\documentclass[12pt]{article}
\usepackage[T1]{fontenc}
\usepackage{lmodern}
\usepackage{latexsym}
\usepackage{authblk}
\usepackage{indentfirst}
\usepackage{amsmath,amssymb}
\usepackage[noadjust]{cite}
\usepackage{tikz}
\usepackage[small]{caption}
\usepackage[hidelinks]{hyperref}
\hypersetup{pdftitle={Clique-cutsets and spanning trees without small-degree stems},
	pdfauthor={Yibo Li}}

\newtheorem{theorem}{Theorem}[section]
\newtheorem{corollary}[theorem]{Corollary}
\newtheorem{lemma}[theorem]{Lemma}
\newtheorem{claim}{Claim}[section]
\newenvironment{proof}{\par\noindent\textbf{Proof.}\enspace}
{\unskip\nobreak\hfill\mbox{\rule{2mm}{2mm}}\par\medskip}
\newcommand{\NC}{\operatorname{NC}}

\begin{document}
	\title{Clique-cutsets and spanning trees without small-degree stems}
\author{Yibo Li\thanks{Email: \href{mailto:liyibo@hue.edu.cn}{\nolinkurl{liyibo@hue.edu.cn}}}}
\affil{School of Mathematics and Statistics, Hubei University of Education, Wuhan 430205, China}
\date{}

	\maketitle
	
	\begin{abstract}
		For an integer $k\ge2$, a $[2,k]$-ST is a spanning tree
		with no vertex of degree $2$ through $k$. Let $G$ be a
		connected graph of sufficiently large order $n$. We prove
		that if $|N_G(x)\cup N_G(y)|\ge(n-2)/2$ for every pair of
		nonadjacent vertices $x,y$, then $G$ has a $[2,k]$-ST
		if and only if it does not belong to either of the two
		exceptional families. In particular, we answer
		Conjecture~1.7 proposed in~\cite{LDL}.

		%One is defined by a vertex cut whose vertices have degrees between $2$ and $k$. The other is defined by the way a clique-cutset is joined to the two components it separates. If $\delta(G)\ge2k-1$, the only obstruction is a clique-cutset of order $2k-1$ whose vertices each have exactly one neighbor outside it. This gives Conjecture~1.7 of Li, Dong and Liu as a special case.
	\end{abstract}
	
	\noindent\textbf{Keywords}\quad $[2,k]$-spanning tree, neighborhood union,
	clique-cutset, minimum degree
	
	\section{Introduction}
	
	All graphs considered in this paper are finite and simple.
	For a graph $G$, let $V(G)$ and $E(G)$ denote its vertex set
	and edge set, respectively. For $v\in V(G)$, let $N_G(v)$
	denote its {\em neighborhood} and let $d_G(v)=|N_G(v)|$.
	We omit the subscript $G$ when the graph is clear, and let
	$\delta(G)$ denote the minimum degree of $G$. For
	$X\subseteq V(G)$, let $N_X(v)=N(v)\cap X$ and
	$d_X(v)=|N_X(v)|$. We write $G[X]$ for the subgraph
	induced by $X$ and $G-X=G[V(G)\setminus X]$. A {\em vertex cut}
	of a connected graph $G$ is a set $Q$ such that $G-Q$ has at least two
	nonempty components. A vertex cut $Q$ of a connected graph $G$ is called
	a {\em clique-cutset} if $G[Q]$ is complete. For a noncomplete graph $G$, we define 
	\[
	\begin{aligned}
		\NC(G)&:=\min\{\,|N(u)\cup N(v)|:
		u,v\in V(G),\ u\ne v,\ uv\notin E(G)\,\},\\
		\sigma(G)&:=\min\{\,d(u)+d(v):
		u,v\in V(G),\ u\ne v,\ uv\notin E(G)\,\}.
	\end{aligned}
	\]

	For an integer $k\ge2$, a spanning tree of $G$ is called a \emph{$[2,k]$-ST}
	if it has no vertex of degree $2$ through $k$. The case $k=2$ is a
	{\em homeomorphically irreducible spanning tree} (HIST). Sufficient conditions for the existence of HISTs have
	been investigated from various perspectives, including
	degree conditions, local connectivity and forbidden
	subgraphs, see~\cite{ABHT90,ItTs22,ChRS12,ChSh13,FuTs20}. For the neighborhood union condition, Li {\em et al.}\cite{LiDH26}
	proved that every connected graph of order $n\ge270$
	with $\NC(G)\ge(n-1)/2$ has a HIST.
	
	Several sufficient conditions for the existence of
	$[2,k]$-STs have been established in terms of minimum
	degree~\cite{FuST25} and degree sum~\cite{FuT25}.
	In particular, Furuya and Tsuchiya proved that, for
	sufficiently large connected graphs, the condition
	$\sigma(G)\ge n-2$ guarantees a $[2,k]$-ST except for one exceptional family. Recently, Li {\em et al.}\cite{LDL} established the following
	neighborhood union condition for the existence of a
	$[2,k]$-ST. 
	
	Let
	$c_k=(\sqrt{k}+\sqrt{2})\sqrt{k(k-1)}$, and let
	$n_1(k)$ be the smallest positive integer $N$ such that
	$n-4c_k\sqrt n-12k+14\ge0$ for every integer $n\ge N$.
	
	\begin{theorem}\cite{LDL}\label{thm:LDL}
		Let $k\ge2$ be an integer, and let $G$ be a connected graph
		of order $n\ge n_1(k)$. If $\delta(G)\ge2k$ and
		\[
		\NC(G)\ge\frac{n-2}{2},
		\]
		then $G$ has a $[2,k]$-ST.
	\end{theorem}
	
	The minimum degree condition $\delta(G)\ge2k$ is best possible. In \cite{LDL}, the authors constructed examples of connected graphs with minimum degree $2k-1$ that have
	no $[2,k]$-ST. In each example, there is a clique-cutset $Q$ of order
	$2k-1$ such that every vertex of $Q$ has exactly one
	neighbor in $V(G)\setminus Q$. This observation led the
	authors to propose Conjecture~1.7 in~\cite{LDL}, which
	states that every sufficiently large connected graph
	$G$ satisfying $\delta(G)\ge2k-1$ and $\NC(G)\ge(n-2)/2$ has a $[2,k]$-ST unless it contains
	such a clique-cutset.

	In this paper, we prove this conjecture and give a more
	general result without a minimum degree condition.
	We first introduce two families of connected graphs,
	$\mathcal F_1$ and $\mathcal F_2$.
	
	For a fixed integer $k\ge2$, let $\mathcal F_1$ be the family of connected graphs having a vertex cut whose vertices all have degrees between $2$ and $k$ in $G$. Such a cut is called a \emph{$k$-blocking set} \cite{FuT25}. Let $\mathcal F_2$ be the family of connected graphs $G\notin\mathcal F_1$ having a clique-cutset $Q$ such that $G-Q$ has exactly two components, say $C_1$ and $C_2$, and every vertex of $Q$ with neighbors in both components has degree at most $k$. For $i=1,2$, let $m_i=\max\{d_{C_i}(v):v\in Q\}$. Without loss of generality, assume that $m_1\le m_2$. We further require that
	\begin{equation}\label{eq:F2}
		|Q|+m_1+m_2\le 2k+1.
	\end{equation}

	\begin{theorem}\label{thm:main}
		Let $k\ge2$ be an integer, and let $G$ be a connected graph of order $n\ge n_1(k)$
		with $\NC(G)\ge(n-2)/2$. Then $G$ has a $[2,k]$-ST if and only if $G\notin\mathcal F_1\cup\mathcal F_2$.
	\end{theorem}
	
	The following corollary confirms Conjecture~1.7 of \cite{LDL} with
	$n_2(k)=n_1(k)$.
	
	\begin{corollary}\label{cor:conjecture}
		Let $k\ge2$ be an integer, and let $G$ be a connected graph of order $n\ge n_1(k)$ with
		$\NC(G)\ge(n-2)/2$ and $\delta(G)\ge2k-1$. Then $G$ has no $[2,k]$-ST
		if and only if it has a clique-cutset $Q$ of order $2k-1$ such that
		every vertex of $Q$ has exactly one neighbor in $V(G)\setminus Q$.
	\end{corollary}

	\section{Preliminaries}\label{sec:pre}
	
	In this section, we first establish a necessary condition for a graph to contain a $[2,k]$-ST.
	
	\begin{lemma}\label{lem:cut}
		Let $G$ be a connected graph, and let $Q$ be a vertex cut of $G$ such that $G-Q$ has exactly two components $C_1$ and $C_2$. Assume that every vertex of $Q$ having neighbors in both $C_1$ and $C_2$ has degree at most $k$. For $i=1,2$, let $m_i=\max_{v\in Q}d_{C_i}(v)$. If $G$ has a $[2,k]$-ST, then
		$|Q|+m_1+m_2\ge 2k+2$.
	\end{lemma}
	
	\begin{proof}
		Suppose that $G$ has a $[2,k]$-ST $T$. Let $P$ be a shortest path in $T$ joining $C_1$ and $C_2$. Let $x$ and $y$ be the first and last vertices of $P$ in $Q$, respectively. If $x=y$, then $x$ has neighbors in both $C_1$ and $C_2$, so $d_G(x)\le k$, contradicting $d_T(x)\ge k+1$. Thus $x\ne y$.
		
		By the choice of $P$, $x$ has no neighbor in $C_2$ in $T$, and $y$ has no neighbor in $C_1$ in $T$. Hence
		\[
		\begin{aligned}
			d_{T[Q]}(x)&\ge d_T(x)-d_{T[C_1]}(x)\ge k+1-m_1,\\
			d_{T[Q]}(y)&\ge d_T(y)-d_{T[C_2]}(y)\ge k+1-m_2.
		\end{aligned}
		\]
		Since $x\ne y$ and $T[Q]$ is a forest, we have
		\[
		\begin{aligned}
			|Q|-1&\ge |E(T[Q])|\ge d_{T[Q]}(x)+d_{T[Q]}(y)-1\\
			&\ge (k+1-m_1)+(k+1-m_2)-1.
		\end{aligned}
		\]
		This yields the desired inequality.
	\end{proof}
	
	For the rest of this section, let $G$ be a connected graph of order $n\ge n_1(k)$ satisfying $\NC(G)\ge(n-2)/2$ and $\delta(G)\le2k-1$. Let $u\in V(G)$ with $d(u)=\delta(G)$. Define
	\[
	S=\{v\in V(G):d(v)<(n-2)/4\},\qquad s=|S|.
	\]
	By the definition of $n_1(k)$, we have
	\[
	n-4c_k\sqrt n-12k+14\ge0.
	\]
	Since $12k-14>0$ for $k\ge2$, it follows that
	\[
	4c_k\sqrt n\le n-12k+14<n.
	\]
	Consequently,
	\begin{equation}\label{eq:order}
		c_k\sqrt n<n/4,\qquad n>16c_k^2\ge64k^2,
	\end{equation}
	where the last inequality follows from $c_k\ge2k$.
	
	\begin{lemma}\label{lem:components}
		The set $S$ is a clique with $1\le |S|\le 2k$. Moreover, every component of $G-S$ has at least $n/2-d(u)$ vertices, and $G-S$ has at most two components.
	\end{lemma}
	
	\begin{proof}
		Suppose that $S$ is not a clique. Then there exist two nonadjacent vertices $x,y\in S$. By the definition of $S$, we have
		\[
		|N(x)\cup N(y)|
		\le d(x)+d(y)
		<\frac{n-2}{2},
		\]
		a contradiction. Hence, $S$ is a clique.
		
		By \eqref{eq:order}, we have  $d(u)\le2k-1<(n-2)/4$, then $u\in S$. Since $S$ is a clique, we obtain
		\[
		1\le s\le d(u)+1\le2k.
		\]
		In particular, $G-S$ is nonempty.
		
		Let $C$ be any component of $G-S$. Since every vertex of $C$ has degree at least $(n-2)/4$ in $G$ and has at most $s$ neighbors in $S$, we have
		\begin{equation}\label{eq:degree-C}
			\delta(G[C])\ge(n-2)/4-s.
		\end{equation}
		By $s\le2k$ and \eqref{eq:order}, we have 
		\[
		\delta(G[C])\ge(n-2)/4-2k>2k-1\ge d(u).
		\]
		Thus, there exists a vertex $w\in C$ such that $uw\notin E(G)$. Otherwise, we have $C\subseteq N(u)$, contradicting $\delta(G[C])>d(u)$.
		
		Since $S$ is a clique and $u\in S$, every neighbor of $w$ in $S$ belongs to $N(u)$. Moreover, all neighbors of $w$ outside $S$ lie in $C$. Hence
		\[
		N(w)\setminus N(u)\subseteq C\setminus\{w\},
		\]
		and consequently
		\[
		\frac{n-2}{2}\le |N(u)\cup N(w)|
		\le d(u)+|C|-1.
		\]
		It follows that
		\[
		|C|\ge\frac n2-d(u)>\frac n3.
		\]
		Since the components of $G-S$ are pairwise disjoint, $G-S$ has at most two components.
	\end{proof}
	
	\vspace{1mm}
	
	We use the following two results in the proof of Lemma~\ref{lem:rooted}.
	
	\vspace{1mm}
	
	\begin{lemma} \cite{FuT25}
		\label{lem:neighbors}
		Let $G$ be a connected graph, $x\in V(G)$, and
		$Y\subseteq V(G)\setminus\{x\}$. For an integer $m\ge0$, if $\delta(G)\ge 2\sqrt{|V(G)|}+m+|Y|$, then there exists a set $L\subseteq N_G(x)\setminus Y$ with
		$|L|=m$ such that $G-L$ is connected.
	\end{lemma}
	
	\vspace{2mm}
	
	For an integer $k\ge2$, let $c_k=(\sqrt{k}+\sqrt{2})\sqrt{k(k-1)}$.
	
	\vspace{2mm}
	
	\begin{lemma}\cite{FuT25}
		\label{lem:tree}
		Let $G$ be a connected graph, and let
		$\emptyset\neq U\subseteq V(G)$. If $\delta(G)\ge c_k\sqrt{|V(G)|}+k|U|-1$, then $G$ has a spanning tree $T$ in which each vertex of $U$ has
		degree at least $k$ and each vertex outside $U$ is a leaf or has
		degree at least $k+1$.
	\end{lemma}
	
	\vspace{1mm}
	
	The following lemma provides a useful construction for the proof of our main theorem.
	
	\vspace{1mm}
	
	\begin{lemma}\label{lem:rooted}
		Let $C$ be a component of $G-S$, and let $R\subseteq V(C)$
		with $1\le |R|\le k+1$. Then $C$ has a spanning tree $T$
		such that
		$d_T(v)\ge k+1$ for every $v\in R$, and $d_T(v)=1$ or $d_T(v)\ge k+1$
		for every $v\in V(C)\setminus R$. Moreover, $C$ has a spanning forest $F$ satisfying the same
		degree conditions, with each component of $F$ containing
		exactly one vertex of $R$.
	\end{lemma}
	\begin{proof}
		Let $X=R\cup N_C(u)$ and $b=|X|$. Write
		$X=\{x_1,\ldots,x_b\}$. Since $S\setminus\{u\}\subseteq N(u)$,
		we have
		\[
		b\le |R|+d(u)-(s-1)\le k+1+2k-1-(s-1)\le3k.
		\]
		By (3), $s\le2k$, and (2), we have
		\begin{equation}\label{eq:4}
			\begin{aligned}
				\delta(G[C])-2\sqrt n
				&\ge \frac{n-2}{4}-2k-2\sqrt n\ge \frac n8-2k-\frac12\ge 3k(k+2)\ge(k+2)b.
			\end{aligned}
		\end{equation}
		Here we use $n>64k^2$ and $b\le3k$.
		
		Set $H_1=C$. Suppose that $H_i$ is connected. By
		\eqref{eq:4} and $i\le b$, we have
		\[
		\begin{aligned}
			\delta(H_i)
			&\ge\delta(G[C])-(i-1)(k+1)\ge2\sqrt n+(k+2)b-(i-1)(k+1)\\
			&\ge2\sqrt{|V(H_i)|}+(k+1)+(b-1),
		\end{aligned}
		\]
		where $H_i=C-\bigcup_{j=1}^{i-1}L_j$.
		Applying Lemma~2.3 with $G=H_i$, $x=x_i$,
		$Y=X\setminus\{x_i\}$ and $m=k+1$, we obtain
		a set $L_i\subseteq N_{H_i}(x_i)\setminus X$ with $|L_i|=k+1$ such that $H_i-L_i$ is connected.
		Set $H_{i+1}=H_i-L_i$. Starting with $H_1=C$
		and repeating this argument for $i=1,\ldots,b$,
		we obtain pairwise disjoint sets $L_1,\ldots,L_b$.
		Consequently, $C-L$ is connected, where
		$L=\bigcup_{i=1}^b L_i$.
		
		Let $B_1,\ldots,B_t$ be the components of
		$C-(L\cup X)$. Contract each $B_i$ into a vertex
		$b_i$, and delete loops and parallel edges. Let $B_1,\ldots,B_t$ be the components of
		$C-(L\cup X)$. Contract each $B_i$ into a vertex
		$b_i$, and delete loops and parallel edges. The resulting graph $H^*$ is connected since
		$C-L$ is connected, with $V(H^*)=X\cup\{b_1,\ldots,b_t\}$.
		Moreover, the vertices $b_1,\ldots,b_t$ are pairwise nonadjacent, and each $b_i$ is adjacent only to vertices
		of $X$.
        
        Clearly, $H^*$ has a spanning tree $T^*$. We also construct a spanning forest $F^*$ of $H^*$ as follows. Start with the vertices of $R$ as isolated vertices and repeatedly add a new vertex together with an edge joining it to the current forest. Since $H^*$ is connected, this yields a
        spanning forest in which each component contains
        exactly one vertex of $R$.

        Now, we have obtained a spanning tree $T^\ast$ and a spanning
        forest $F^\ast$ of $H^\ast$. We next expand these
        structures to subgraphs of $C-L$. For convenience, let
        $K^\ast$ denote either $T^\ast$ or $F^\ast$, and let $K$
        denote its corresponding expansion in $C-L$.
        For each edge $b_ix\in E(K^\ast)$ with $x\in X$,
        choose a corresponding edge between $B_i$ and $x$
        in $G$. Keep the edges of $K^\ast$ whose endpoints
        both belong to $X$. Let $K$ be the subgraph of $C-L$
        formed by these edges.
        
        For each component $B_i$ of $C-(L\cup X)$, let
        \[
        U_{B_i}(K)=\{v\in V(B_i):v\text{ is incident with an edge of }K\}.
        \]
        Since $K^*$ is a spanning tree or a spanning forest
        of $H^*$, each $b_i$ has at least one neighbor in $X$. Moreover, the vertices $b_1,\ldots,b_t$
        are pairwise nonadjacent, and each $b_i$ is adjacent
        only to vertices of $X$. Since $H^*$ is simple, each vertex of $X$ is incident with at most one
        selected edge joining it to $B_i$. Consequently,
        $1\le|U_{B_i}(K)|\le b$.

We verify that Lemma~2.4 applies to $G[B]$ and
$U_B(K)$. For every $v\in V(B)$, we have $v\in C-X$
and $uv\notin E(G)$. Since $S$ is a clique and $u\in S$,
every neighbor of $v$ in $S$ belongs to $N(u)$.
Furthermore, every neighbor of $v$ outside $C$ belongs
to $S$. Hence $N(v)\setminus N(u)\subseteq V(C)$. By the neighborhood union condition, we obtain
\begin{equation}\label{eq:outside}
d_C(v)\ge |N(u)\cup N(v)|-d(u)
\ge \frac{n-2}{2}-d(u)
\ge \frac n2-2k.
\end{equation}

Since every neighbor of $v$ in $C\setminus B$ belongs
to $L\cup X$, we have
\[
\begin{aligned}
\delta(G[B])
&\ge \frac n2-2k-(k+2)b\ge c_k\sqrt n+kb\ge c_k\sqrt{|V(B)|}+k|U_B(K)|-1.
\end{aligned}
\]
The second inequality follows from
$c_k\sqrt n<n/4$ and
$n/4-2k\ge(2k+2)b$, which follow from
\eqref{eq:order} and $b\le3k$. The last inequality
follows from $|V(B)|\le n$ and $|U_B(K)|\le b$.

By Lemma~2.4, each $B_i$ has a spanning tree
$T_{B_i}(K)$ such that $d_{T_{B_i}(K)}(v)\ge k$
for every $v\in U_{B_i}(K)$, and every
$v\in V(B_i)\setminus U_{B_i}(K)$ has degree $1$
or at least $k+1$ in $T_{B_i}(K)$. Replacing each contracted vertex $b_i$ in $K^\ast$
by $T_{B_i}(K)$, and retaining all edges in $E(K)$,
gives a spanning tree of $C-L$ when $K^\ast=T^\ast$,
and a spanning forest of $C-L$ when $K^\ast=F^\ast$.
Every vertex of $U_{B_i}(K)$ has at least one selected
edge to $X$, and hence has degree at least $k+1$.
The degrees of vertices in
$V(B_i)\setminus U_{B_i}(K)$ remain unchanged.

Finally, attach each vertex of $L_i$ as a leaf to
$x_i$, for $i=1,\ldots,b$. Since $|L_i|=k+1$ and
the sets $L_1,\ldots,L_b$ are pairwise disjoint and
disjoint from $X$, every vertex of $X$ has degree
at least $k+1$. Thus, the resulting spanning tree and spanning forest
satisfy the required degree conditions. Moreover,
each component of the forest contains exactly one
vertex of $R$. This completes the proof.
	\end{proof}

	\section{Proof of \texorpdfstring{Theorem~\ref{thm:main}}{the main theorem}}
	In this section, we will prove Theorem~\ref{thm:main}.
	Throughout this section, we always assume that $G$ is a graph satisfying the
	hypotheses of Theorem~\ref{thm:main}.
	
	Suppose first that $G$ has a $[2,k]$-ST $T$. If
	$G\in\mathcal F_1$, then $G$ has a $k$-blocking set $B$.
	Since every vertex of $B$ has degree between $2$ and $k$
	in $G$, each vertex of $B$ must be a leaf of $T$. Thus,
	$T-B$ is connected, contradicting the fact that $B$
	is a vertex cut. If $G\in\mathcal F_2$, then
	Lemma~\ref{lem:cut} implies that
	$|Q|+m_1+m_2\ge 2k+2$, contrary to \eqref{eq:F2}.
	Hence, $G\notin\mathcal F_1\cup\mathcal F_2$.
	
	For the converse, suppose that $G$ has no $[2,k]$-ST. Let $u\in V(G)$ with $d(u)=\delta(G)$. By Theorem~\ref{thm:LDL}, we have $d(u)\le2k-1$. Recall that $S=\{v\in V(G):d(v)<(n-2)/4\}$ and $s=|S|$. By
	Lemma~\ref{lem:components}, $S$ is a clique and $G-S$ has
	at most two components.

	We first establish the following claim.
	
	\begin{claim}\label{cl:one}
		If $z\in S$ is adjacent to a vertex in every component of
		$G-S$, then $d(z)\le k$.
	\end{claim}
	\begin{proof}
		Suppose that $d(z)\ge k+1$. Since $S$ is a clique,
		$z$ is adjacent to every vertex of $S\setminus\{z\}$.
		Let $C_1,\ldots,C_q$ be the components of $G-S$,
		where $q\in\{1,2\}$. Let
		$r=\max\{q,k+2-s\}$. Since $N_{C_i}(z)\neq\emptyset$ for each $i\in[q]$ and
		$d(z)\ge k+1$, there exists $R\subseteq N(z)\setminus S$ such that $|R|=r$ and $R\cap V(C_i)\ne\emptyset$ for $i\in[q]$.
		Then $1\le |R|\le k+1$ and $s-1+|R|\ge k+1$.
		
		By Lemma~\ref{lem:rooted}, each $C_i$ has a spanning
		forest $F_i$ such that
		\[
		d_{F_i}(v)\ge k+1\quad(v\in R\cap V(C_i)),
		\]
		and
		\[
		d_{F_i}(v)=1\quad\text{or}\quad d_{F_i}(v)\ge k+1
		\quad(v\in V(C_i)\setminus R),
		\]
		with each component of $F_i$ containing exactly one
		vertex of $R$.
		
		Let
		\[
		T=\{zv:v\in S\setminus\{z\}\}
		\cup\bigcup_{i=1}^{q}F_i
		\cup\{zv:v\in R\}.
		\]
		Since each component of $F_i$ contains exactly one
		vertex of $R$, $T$ is a spanning tree of $G$.
		Moreover,
		\[
		d_T(z)=s-1+|R|\ge k+1,
		\]
		\[
		d_T(v)=1\quad(v\in S\setminus\{z\}),
		\]
		and
		\[
		d_T(v)=1\quad\text{or}\quad d_T(v)\ge k+1
		\quad(v\in V(G)\setminus S).
		\]
		Thus, $T$ is a $[2,k]$-ST, a contradiction.
	\end{proof}
	
	We consider the following two cases.

	\medskip
	
	\noindent\textbf{Case 1.} \textit{$G-S$ is connected.}
	
	\medskip
	
	In this case, let $P=\{x\in S:N_{G-S}(x)\ne\emptyset\}$. Since $G$ is connected, $P\ne\emptyset$. By
	Claim~\ref{cl:one}, $d(x)\le k$ for every $x\in P$.
	
	If $P\ne S$, then $P$ is a vertex cut separating
	$S\setminus P$ from $G-S$. Since $S$ is a clique, every $x\in P$ has a neighbor in $S\setminus P$.
	Moreover, by the definition of $P$, $d_{G-S}(x)\ge1$. Hence, $2\le d(x)\le k$ for each $x\in P$. Thus, $P$ is a $k$-blocking set, contrary to $G\notin\mathcal F_1$.
	
	If $P=S$, then $d(x)\ge s$ for every $x\in S$, and hence
	$s\le k$. For each $x\in S$, choose
	$w_x\in N_{G-S}(x)$, and let $R=\{w_x:x\in S\}$. By Lemma~\ref{lem:rooted}, $G-S$ has a spanning tree $T_0$ such that $d_{T_0}(v)\ge k+1$ for each $v\in R$, $d_{T_0}(v)=1$ or $d_{T_0}(v)\ge k+1$ for each $v\in V(G-S)\setminus R$.
	Let $T=T_0\cup\{xw_x:x\in S\}$. Then $T$ is a spanning tree of $G$, and $d_T(x)=1$ for each $x\in S$, $d_T(v)=1$ or $d_T(v)\ge k+1$ for each $v\in V(G-S)$. Thus, $T$ is a $[2,k]$-ST, a contradiction.
	
	\medskip
	
	\noindent\textbf{Case 2.} \textit{$G-S$ is disconnected.}
	
	\medskip
	
	In this case, let $C_1$ and $C_2$ be two components of $G-S$.
	For $i=1,2$, let $P_i=\{v\in S:N_{C_i}(v)\ne\emptyset\}$ and
	$m_i=\max_{x\in S}d_{C_i}(x)$.
	Clearly, $P_i$ is a vertex cut separating $C_i$ and $C_{3-i}$. Moreover, we have $d(x)\ge2$ for every $x\in P_i$. If $d(x)\le k$ for all $x\in P_i$, then $P_i$ is a
	$k$-blocking set, contrary to $G\notin\mathcal F_1$.
	Thus, for each $i\in\{1,2\}$, there exists
	$x_i'\in P_i$ such that $d(x_i')\ge k+1$. By Claim~\ref{cl:one}, $N_{C_{3-i}}(x_i')=\emptyset$.
	Since $S$ is a clique, we have
	\[
	d_{C_i}(x_i')=d(x_i')-(s-1)\ge k+2-s.
	\]
	On the other hand, for any $y\in S$ with $d(y)\le k$,
	\[
	d_{C_i}(y)\le d(y)-(s-1)\le k+1-s.
	\]
	Thus, there exists $x_i\in P_i$ such that
	$d_{C_i}(x_i)=m_i$ and $d(x_i)\ge k+1$. By
	Claim~\ref{cl:one}, $N_{C_{3-i}}(x_i)=\emptyset$.
	Hence, $x_1\ne x_2$ and
	\begin{equation}\label{eq:individual}
		s-1+m_i=d(x_i)\ge k+1,\qquad i=1,2.
	\end{equation}
	
	Suppose that $s+m_1+m_2\ge2k+2$. For $i=1,2$, let $p_i=\min\{k,m_i\}$, and let $R_i\subseteq N_{C_i}(x_i)$ with
	$|R_i|=p_i$.
	By \eqref{eq:individual} and the above inequality,
	\[
	(k-p_1)+(k-p_2)\le s-2.
	\]
	Hence, there is a partition
	$S\setminus\{x_1,x_2\}=A_1\cup A_2$ such that
	$|A_i|\ge k-p_i$, $i=1,2$.
		
	By Lemma~\ref{lem:rooted}, for each $i=1,2$,
	$C_i$ has a spanning forest $F_i$ in which $d_{F_i}(v)\ge k+1$ for $v\in R_i$,
	and $d_{F_i}(v)=1$ or $d_{F_i}(v)\ge k+1$ for $v\in V(C_i)\setminus R_i$ with each component containing exactly one vertex of $R_i$.
	
	Let $T=F_1\cup F_2
	\cup\{x_1x_2\}
	\cup\{x_iv:v\in A_i\cup R_i,\ i=1,2\}$. Then $T$ is a spanning tree of $G$. Moreover,
	\[
	d_T(x_i)=1+|A_i|+p_i\ge k+1,
	\qquad i=1,2,
	\]
	and
	\[
	d_T(v)=1\quad(v\in S\setminus\{x_1,x_2\}).
	\]
	The remaining vertices have degree $1$ or at least $k+1$.
	Thus, $T$ is a $[2,k]$-ST, a contradiction. Consequently,
	\begin{equation}\label{eq:cut}
		s+m_1+m_2\le2k+1.
	\end{equation}
		Combining \eqref{eq:individual} and \eqref{eq:cut},
		we obtain $3\le s\le2k-1$.
	
	By Claim~\ref{cl:one}, every vertex of $S$ having
	neighbors in both $C_1$ and $C_2$ has degree at most $k$.
	Together with \eqref{eq:cut}, this shows that $Q=S$
	satisfies the definition of $\mathcal F_2$.
	This contradicts $G\notin\mathcal F_2$.
	
	\section{Further Results}
	
	\noindent\textbf{Proof of Corollary~\ref{cor:conjecture}.}
	Since $\delta(G)\ge2k-1>k$, Theorem~\ref{thm:main}
	implies that $G\in\mathcal F_2$ if $G$ has no
	$[2,k]$-ST. Let $Q$ be the corresponding clique-cutset,
	and let $C_1,C_2$ be the two components of $G-Q$.
	By the definition of $\mathcal F_2$, no vertex of $Q$
	has neighbors in both $C_1$ and $C_2$. For $i=1,2$,
	let $m_i=\max_{x\in Q}d_{C_i}(x)$, and choose
	$x_i\in Q$ with $d_{C_i}(x_i)=m_i$. Then
	$|Q|+m_i=d(x_i)+1\ge2k$, for $i=1,2$.
	Together with $m_1,m_2\ge1$ and \eqref{eq:F2},
	this gives $|Q|=2k-1$ and $m_1=m_2=1$.
	Hence, every vertex of $Q$ has exactly one neighbor
	in $V(G)\setminus Q$.
	
	Conversely, suppose that \(G\) has a clique-cutset \(Q\) of order \(2k-1\) such that every vertex of \(Q\) has exactly one neighbor in \(V(G)\setminus Q\). Suppose, to the contrary, that \(G\) has a \([2,k]\)-ST \(T\). Choose two distinct components of \(G-Q\), and let \(P\) be a shortest path in \(T\) joining them. Let \(x\) and \(y\) be the first and last vertices of \(Q\) on \(P\), respectively. Since every vertex of \(Q\) has exactly one neighbor outside \(Q\), we have \(x\ne y\). Moreover, since \(T\) is a \([2,k]\)-ST and each vertex of \(Q\) has at most one neighbor outside \(Q\), we have $d_{T[Q]}(x)\ge k$ and $d_{T[Q]}(y)\ge k$. Since \(T[Q]\) is a forest, it follows that
	\[
	2k-2=|Q|-1\ge |E(T[Q])|
	\ge d_{T[Q]}(x)+d_{T[Q]}(y)-1
	\ge 2k-1,
	\]
	a contradiction. Thus, \(G\) has no \([2,k]\)-ST.
	\hfill\mbox{\rule{2mm}{2mm}}
	
	\begin{corollary}\label{cor:lower}
		Let $k\ge3$ be an integer, and let $G$ be a connected graph
		of order $n\ge n_1(k)$ with
		$\NC(G)\ge(n-2)/2$ and $\delta(G)\ge2k-2$.
		Then $G$ has no $[2,k]$-ST if and only if $G$ has a
		clique-cutset $Q$ such that $G-Q$ has exactly two components,
		no vertex of $Q$ has neighbors in both components, and
		\[
		(|Q|,m_1,m_2)\in
		\{(2k-1,1,1),(2k-2,1,1),(2k-2,1,2),(2k-3,2,2)\},
		\]
		where $m_1\le m_2$ are as defined in $\mathcal F_2$.
	\end{corollary}
	\begin{proof}
		Since $\delta(G)\ge2k-2>k$, Theorem~\ref{thm:main}
		implies that $G\in\mathcal F_2$ if it has no $[2,k]$-ST.
		Thus, $G-Q$ has exactly two components, no vertex of $Q$
		has neighbors in both components, and \eqref{eq:F2} holds.
		Moreover, $\delta(G)\ge2k-2$ gives
		$|Q|+m_1\ge2k-1$. Thus, $1\le m_1\le m_2\le2$.
		If $m_2=1$, then $m_1=1$ and
		$2k-2\le|Q|\le2k-1$. If $m_2=2$, then
		$|Q|+m_1\ge2k-1$ and \eqref{eq:F2} imply
		$|Q|+m_1=2k-1$. Hence, the four possible triples are
		\[
		(|Q|,m_1,m_2)\in
		\{(2k-1,1,1),(2k-2,1,1),
		(2k-2,1,2),(2k-3,2,2)\}.
		\]
		
		Conversely, each of these four triples satisfies
		\eqref{eq:F2} and the corresponding lower bound
		$|Q|+m_1\ge2k-1$. By the definition of $\mathcal F_2$
		and Theorem~\ref{thm:main}, $G$ has no $[2,k]$-ST.
	\end{proof}
	
	The four possibilities in Corollary~\ref{cor:lower} are sharp,
	as shown by the following construction.
	Let $Q$, $C_1$ and $C_2$ be three disjoint clique, where
	$|Q|$ is as specified and $|C_1|=|C_2|$ are sufficiently
	large. If $|Q|=2k-1$, choose one vertex $x\in Q$ and partition $Q\setminus\{x\}$ into two nonempty sets $Q_1$ and $Q_2$. Join each vertex of $Q_1$
	to $m_1$ vertices of $C_1$, and each vertex of $Q_2$
	to $m_2$ vertices of $C_2$. Add no other edges.
	
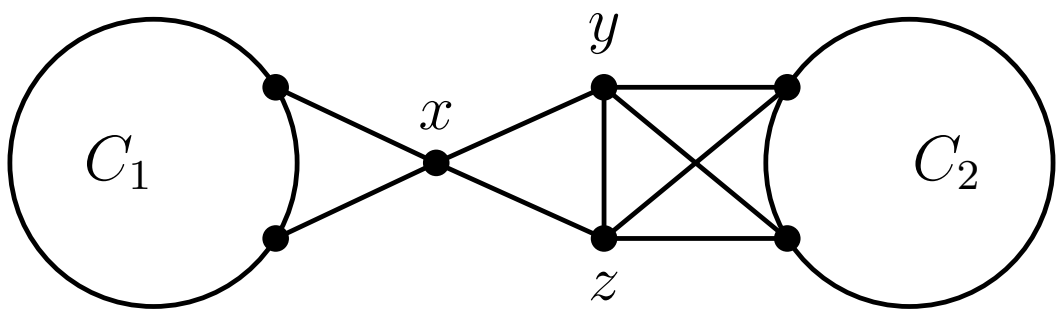
\begin{figure}[htbp]
		\centering
		\begin{tikzpicture}[
			line width=.9pt,
			vertex/.style={circle,fill=black,inner sep=0pt,minimum size=5pt}
			]
			\draw (-2.5,0) circle (.95);
			\draw (2.5,0) circle (.95);
			
			\node at (-2.73,0) {$C_1$};
			\node at (2.75,0) {$C_2$};
			
			\coordinate (a) at (-1.692,.50);
			\coordinate (c) at (-1.692,-.50);
			\coordinate (b) at (1.692,.50);
			\coordinate (d) at (1.692,-.50);
			\coordinate (x) at (-.63,0);
			\coordinate (y) at (.48,.50);
			\coordinate (z) at (.48,-.50);
			
			\draw
			(a)--(x)--(c)
			(x)--(y)--(z)--cycle
			(y)--(b) (y)--(d)
			(z)--(b) (z)--(d);
			
			\foreach \v in {a,b,c,d,x,y,z}
			\node[vertex] at (\v) {};
			
			\node[above=2pt] at (x) {$x$};
			\node[above=2pt] at (y) {$y$};
			\node[below=2pt] at (z) {$z$};
		\end{tikzpicture}
		\caption{A construction with $k=3$ and $\delta(G)=4$.}
		\label{fig:example}
	\end{figure}

	By the choice of $(|Q|,m_1,m_2)$, we have
	$\delta(G)=2k-2$ and $\max_{x\in Q}d_{C_i}(x)=m_i$ for $i=1,2$. Moreover, every nonadjacent pair $x,y$ satisfies
	\[
	|N(x)\cup N(y)|\ge
	\begin{cases}
		|Q|+|C_i|-2, & x\in Q,\ y\in C_i,\\
		|C_1|+|C_2|-2, & x\in C_1,\ y\in C_2.
	\end{cases}
	\]
	Since $|C_1|=|C_2|$ are sufficiently large, both
	bounds are at least $(n-2)/2$. Finally,
	$|Q|+m_1+m_2\le2k+1$, so Lemma~\ref{lem:cut}
	implies that $G$ has no $[2,k]$-ST.
	Figure~\ref{fig:example} illustrates the construction
	for $k=3$ and $(|Q|,m_1,m_2)=(3,2,2)$.

	\vspace{3mm}
	
	For $k=2$, Theorem~\ref{thm:main} gives the following characterization
	at $n_1(2)=276$.

	\begin{corollary}\label{cor:hist}
		Let $G$ be a connected graph of order $n\ge276$ with
		$\NC(G)\ge(n-2)/2$. Then $G$ has a HIST if and only if
		it has neither a vertex cut consisting of vertices of
		degree two nor a triangle whose deletion disconnects $G$
		and each of whose vertices has at most one neighbor in
		$V(G)\setminus Q$.
	\end{corollary}
	
	\begin{proof}
		Suppose that $G$ has no HIST and no vertex cut consisting
		of vertices of degree two. By Case~2, we have
		$3\le |Q|\le2k-1=3$, where $k=2$. Thus, $Q$ is a
		triangle. By \eqref{eq:cut} and $m_1,m_2\ge1$, we obtain
		$m_1=m_2=1$. Hence, no vertex of $Q$ has more than one
		neighbor in $V(G)\setminus Q$.
		
		Conversely, suppose that $G$ has a vertex cut $S$
		consisting of vertices of degree two. If $G$ has a
		HIST $T$, then every vertex of $S$ must be a leaf
		of $T$. Thus, $T-S$ is connected, contradicting
		that $S$ is a vertex cut. Next we suppose that $G$ has a triangle $Q$ whose deletion
		disconnects $G$ and each vertex of $Q$ has at most
		one neighbor in $V(G)\setminus Q$. Suppose that $G$
		has a HIST $T$. Choose two distinct components of
		$G-Q$, and let $P$ be a shortest path in $T$ joining
		them. Let $x$ and $y$ be the first and last vertices
		of $Q$ on $P$, respectively. Since each vertex of $Q$
		has at most one neighbor outside $Q$, we have $x\ne y$.
		Moreover, since $T$ is a HIST, we have
		$d_{T[Q]}(x),d_{T[Q]}(y)\ge2$.
		Since $T[Q]$ is a forest, this gives
		\[
		\begin{aligned}
			2=|Q|-1
			&\ge |E(T[Q])|\ge d_{T[Q]}(x)+d_{T[Q]}(y)-1\ge 3,
		\end{aligned}
		\]
		a contradiction.
	\end{proof}

	\section*{Declaration of Competing Interest}
	The author declares no known competing financial interests or
	personal relationships that could have appeared to influence the work
	reported in this paper.
	
	\end{document}